\documentclass[a4paper,12pt,ceqn,leqno]{article}
\usepackage{t1enc,amssymb,amsbsy,amsthm,amsmath,amsfonts,pgf,dsfont,paralist}

\newtheorem{theo}{Theorem}
\newtheorem{prop}{Proposition}

\newtheorem{cor}{Corollary}
\newtheorem{lm}{Lemma}
\newtheorem{nota}{Notation}

\newtheorem*{rk*}{Remark}

\newtheorem{df}{Definition}

\newenvironment{prf}{\textbf{Proof}}{\flushright$\blacksquare$\\}

\newcommand{\beq}{\begin{equation*}}
\newcommand{\eeq}{\end{equation*}}
\newcommand{\dis}{\displaystyle}

\begin{document}
\author{Marc Malric}

\vskip3cm

\title{\textbf{\Large Density of paths of iterated Lévy transforms of Brownian
motion}} \maketitle

\noindent \textbf{\underline{Abstract :}} The Lévy transform of a
Brownian motion B is the Brownian motion $\displaystyle
B'_t=\int_0^tsgn \,(B_s)\,dB_s$. Call $T$ the corresponding
transformation on the Wiener space $W$. We establish that a. s. the
orbit of $w(\in W)$ under $T$ is dense in $W$ for the compact
uniform convergence topology.

\newpage
\section{Lévy raisings, B-raised Brownian motions and related tools.}
Let $(\Omega, \mathcal{F}, \mathbf{P})$ be the probability space
where all random elements are defined, and $(\mathbf{W},
\mathcal{W}, \pi)$ the Wiener space. Any measurable map from
$\Omega$ to any measure space, defined $\mathbf{P}-$a.e., will be
called a random variable. If $X$ is a r.v. with values in some
measurable space $A$, the probability measure $\mathbf{P} \circ
X^{-1}$ on $A$ is called the law of $X$, and denoted by
$\mathcal{L}(X)$. For instance, a $\mathbf{W}-$valued r.v. with law
$\pi$ is a Brownian motion.\\
The Lévy transform $\mathbf{T} : \mathbf{W} \longrightarrow
\mathbf{W}$ is defined $\pi-$a.e. and preserves $\pi$. Given a
Brownian motion $B$, we denote by $B^n$ its $n-$Lévy iterate, that
is, the Brownian motion $B^n = \mathbf{T}^n \circ B$.\\
\\
From now on, $T>0$, $\varphi \in \mathcal{C}([0,T], \mathbb{R})$ and $\varepsilon > 0$ are
fixed, and $B$ a Brownian motion. The goal is to prove that
the event $E = E^{\varepsilon} = \{\forall n \geq 0 , \ \|B^n - \varphi \|_{\infty} >
\varepsilon \}$ is negligible, where $\|f\|_{\infty} = \|f_{| [0,T]}\|_{\infty }$. It suffices in fact to show
$\mathbf{P}(E) < \varepsilon$, because $E^{\varepsilon_1} \subset
E^{\varepsilon_2}$ when $\varepsilon_2 < \varepsilon_1$.\\
\\
The idea is to construct from B another stochastic process $\Gamma :
\Omega \longrightarrow \mathbf{W}$, which depends on $T$, $\varphi$ and
$\varepsilon$, and has the following three properties :
\begin{enumerate}
\item \label{l1} The law of the process $\Gamma$, i.e., the probability $F
\longrightarrow \mathbf{P}(\Gamma^{-1}F)$ on $(\mathbf{W},
\mathcal{W})$, is absolutely continuous w.r.t. the law $\pi$ of $B$.
\item \label{l2} For some deterministic $r \geq 0$, one has $\Gamma^r = B^r$, that
is, $\mathbf{T}^r \circ \Gamma = \mathbf{T}^r \circ B$.
\item \label{l3} $\mathbf{P} \left(\forall n \geq 0 \|\Gamma^n - \varphi \|_{\infty} > \varepsilon \right)
<\varepsilon$.
\end{enumerate}
Property (i) implies that $\mathbf{T} \circ \Gamma$ can be (almost
everywhere) defined, in spite of $\mathbf{T}$ not being everywhere
defined. Indeed, if $\mathbf{T'} : \mathbf{W} \longrightarrow
\mathbf{W}$ is another version of $\mathbf{T}$, that is if
$\mathbf{T'} = \mathbf{T}$ a.e., the set $\{\mathbf{T'} \neq
\mathbf{T}\}$ is $\pi-$negligible; hence, by (i),
$\Gamma^{-1}\{\mathbf{T'} \neq \mathbf{T}\}$ is
$\mathbf{P}-$negligible, and $\mathbf{T'} \circ \Gamma = \mathbf{T}
\circ \Gamma \mathbf{P}-$a.s.. Similarly, one can define the
stochastic processes $\Gamma^n = \mathbf{T}^n \circ \Gamma$, which
verify $\Gamma^0 = \Gamma$ and $\mathbf{T} \circ \Gamma^n =
\Gamma^{n+1}$.\\

\begin{prop}
For fixed $G \in \mathcal{W}$ and $\varepsilon > 0$, let us suppose
that there exists a stochastic process $\Gamma : \Omega
\longrightarrow \mathbf{W}$ satisfying properties (i), (ii) and
\begin{equation}
\mathbf{P} \left(\forall n \geq 0, \Gamma^n \in G \right) <
\varepsilon.\end{equation} Then :
\begin{equation}
\mathbf{P} \left( \forall n \geq 0, B^n \in G \right) <\varepsilon.
\end{equation}
\end{prop}
~\\
\begin{prf}
Take $G \in \mathcal{W}$ and put $F = \bigcap_{n \geq 0}
\mathbf{T}^{-n}G$. Then, for $r \geq 0$,
\begin{equation}
\mathbf{T}^{-r}F = \bigcap_{n\geq r} \mathbf{T}^{-n}G \supset
\bigcap_{n \geq 0} \mathbf{T}^{-n}G = F.
\end{equation}
But these two sets, $F$ and $\mathbf{T}^{-r}F$, included in one
another, have the same $\pi-$probability by $\mathbf{T}-$invariance;
so equality $F = \mathbf{T}^{-r}F$ holds up to $\pi-$negligibility.
As the laws of $\Gamma$ and $B$ are absolutely continuous w.r.t.
$\pi$ (this is where (i) is used), we have $\Gamma^{-1}(F) =
\Gamma^{-1} (\mathbf{T}^{-r}F)$ and $B^{-1}F =
B^{-1}(\mathbf{T}^{-r}F)$ up to $\mathbf{P}-$negligible events. In
other words, almost surely, we have $\{\Gamma \in F\} = \{\Gamma^r
\in F\}$ and $\{ B \in F\} = \{B^r \in F\}$. Consequently, choosing
$r$ given by (ii) and using $\Gamma^r = B^r$, we have $\{\Gamma \in
F \} = \{B \in F\}$ a.s.. That is to say :
\begin{equation}
\mathbf{P} \left( \forall n \geq 0, \Gamma^n \in G \right) =
\mathbf{P} \left( \forall n \geq 0, B^n \in G\right)
\end{equation}
\end{prf}~
\\
Specializing $G = \{w \in \mathbf{W}$ ; $\|w - \varphi \|_{\infty} >
\varepsilon\}$, we obtain : \begin{equation*} \mathbf{P}(E) <
\varepsilon
\end{equation*}
\\
Proposition 1 reduces the proof of the approximation
theorem to the construction of a process $\Gamma$ verifying (i),
(ii), and (iii). We shall first choose $r$ in a suitable way, then
work backwards, in $r$ steps, from $\Gamma^r = B^r$ to $\Gamma =
\Gamma^0$; each step (called a \textit{Lévy raise}) will construct
$\Gamma^{n-1}$ from its Lévy transform $\Gamma^n = \mathbf{T} \circ
\Gamma^{n-1}$. The sequence $(\Gamma^r, \Gamma^{r-1}, ...,
\Gamma^0)$ is given a name :

\begin{df}
Given $r \in \mathbb{N}$, a sequence $(\Gamma^r , \Gamma^{r-1} ,...,
\Gamma^0)$ is called a sequence of $B-$raised Brownian motions of
index $r$ if each $\Gamma^n$ is a $\mathbf{W}-$valued r.v. with law absolutely continuous w.r.t.
$\pi$, if $\Gamma^r = B^r$, and if we have $\Gamma^n = \mathbf{T}
\circ \Gamma^{n-1}$ for $0 < n \leq r$.
\end{df}
\textit{In fact, for convenience of exposition, \textit{let us enlarge the
filtered probability space $\Omega$}, we suppose it contains the
whole sequence
$(B'^n)_{n \in \mathbb{N}}$ of the Brownian iterates of $B'$}, B.M. independent from $B$.\\
\\
So we can assert :\\
\begin{cor}
To prove the approximation theorem, it suffices to exhibit
a sequence $(\Gamma^r, \Gamma^{r-1}, ..., \Gamma^0)$ of $B-$raised
Brownian motions of index $r$ such that
\begin{equation}
\mathbf{P} \left( \|\Gamma^n - \varphi \|_{\infty} < \varepsilon \quad for \ some \
n \ \in \{0, ..., r\} \right) > 1 - \varepsilon.
\end{equation}
\end{cor}
~\\
\begin{prf}
Properties (i) and (ii) of Proposition 1 are granted by the
definition of a sequence of $B-$raised Brownian motions, and (iii)
is implied by (1).
\end{prf}~
\\
A Lévy raise starts with a given $\mathbf{W}-$valued r.v.
$\Gamma^n$, and yields some r.v. $\Gamma^{n-1}$ with Lévy transform
$\Gamma^n$. Given a $\mathbf{W}-$valued r.v. $V$, how can one find a
r.v. $U$ such that $V = \mathbf{T} \circ U$? Knowing $V$ is
equivalent to knowing $|U|$, so to define $U$ one only needs to
decide which sign is assigned to each excursion of $|U|$ away from
zero. To make this rigorous, we need a formal definition of the
excursions of a path and of their signs.

\begin{nota}
For $w \in \mathbf{W}$ and $q > 0$, denote by $Z(w) = \left\{s \geq
0 / w(s) = 0 \right\}$ the set of zeros of $w$, and define $g_q(w) =
\sup \left([0,q] \cap Z(w)\right) \geq 0$ (last zero before q) and
$d_q(w) = \inf \left( [q,\infty] \cap Z(w) \right) \neq \infty$
(first zero after q).
\end{nota}

Fix a dense sequence $(q_n)$ in $[0, \infty]$. To each $w \in
\mathbf{W}$, we can attach the sequence $(e_p)$ of disjoint, open
intervals obtained from the sequence

\begin{equation}
\left((g_{q_1}, d_{q_1}), (g_{q_2}, d_{q_2}), ..., (g_{q_n},
d_{q_n}), \dots \right)
\end{equation}

by deleting an interval whenever it already occurs earlier in the
sequence. The $e_p$ are the \textit{excursion intervals} of $w.
\pi-$almost surely, there are infinitely many of them, and they are
the connected components of the open set $[0, \infty] \backslash
Z(w)$. The interval $e_p(w)$ will be called the \textit{p-th
excursion interval} of $w$; $e_p$ is an interval-valued measurable
map, defined on $(\mathbf{W}, \mathcal{W})$ up to
$\pi-$negligibility.\\
Since $w$ does not vanish on $e_p(w)$, its sign is constant on this
interval; this sign will be denoted by $\mathbf{S}_p(w)$, and the
sequence $(\mathbf{S}_p)$ will be called $\mathbf{S}$. If $B$ is a
Brownian motion, the sequence of r.v. $\mathbf{S} \circ B =
(\mathbf{S}_p \circ B)$ is a \textit{coin-tossing}; this means, it
is an i.i.d. sequence, with each r.v. $\mathbf{S}_p \circ B$
uniformly distributed on the set $\{-1, +1\}$. Moreover,
$\mathbf{S}_p \circ B$ and $|B|$ are independent. (See Chap. XII of
[R,Y]).

\begin{lm}
Define $\mathbf{I} : \mathbf{W} \longrightarrow \mathbf{W}$  by
$\mathbf{I}(w)(s) = \inf_{[0, s]}w$ ; that is, $|w| = \mathbf{T}w -
\mathbf{I}\mathbf{T}w$ for $\pi-$a.a. $w$.
\end{lm}
~\\
\begin{prf}
Fix $s \geq 0$. On $[0, s]$, $B^1 = |B| - L \geq - L_s = L_{g_s} =
B_{g_s}^1$. So $B_{g_s}^1 = \inf_{[0, s]} B^1$, and $|B_s| = B_s^1 +
L_s = B_s^1 - \inf_{[0, s]} B^1$.
\end{prf}~
\\
\begin{nota}
If $f$ is a continuous function on $[u,v]$, we call $\arg\min_{[u,v]} f$, resp. $\arg\max_{[u,v]} f$, the largest $t \in [u,v]$ such that $f(t) = \min_[u,v] f$, resp. $f(t) =\max_[u,v] f$.
\end{nota}~
\\
\begin{lm}
Let $A$, $A'$ and $A''$ be three measure spaces; let $\mu_1$ and
$\mu_2$ be two measures on A, $f$ a measurable map from $A$ to $A'$,
and $\nu$ a measure on $A''$. If $\mu_1 \ll \mu_2$, then
\begin{enumerate}
\item \label{l1} $\mu_1 \circ f^{-1} \ll \mu_2 \circ f^{-1}$;
\item \label{l2} $ \mu_1 \otimes \nu  \ll \mu_2 \otimes \nu$.
\end{enumerate}
\end{lm}

~\\
\begin{prf}
(i) If $F \subset A'$ is measurable and if $(\mu_2 \circ f^{-1})(F)
= 0$, then $\mu_2(f^{-1} F) = 0$, so $\mu_1 \circ f^{-1}(F) =
\mu_1(f^{-1} F) = 0$.
\\
(ii) If a measurable subset $F$ of $A \times A''$ is negligible for
$\mu_2 \otimes \nu$, then $\nu-$almost all its sections $F_y$ verify
$\mu_2(F_y) = 0$. Hence they also verify $\mu_1(F_y) = 0$, and
consequently $(\mu_1 \otimes \nu)(F) = \int \mu_1(F_y)\nu(dy) = 0$.
\end{prf}~
\\

\begin{lm}
Let $\tau = (\tau_p)$ be coin-tossing, $\tau' = (\tau'_p)$ a r.v.
with values in $\{-1 , 1\}^{\mathbb{N}}$ such that $\tau'_p =\tau_p$
for all but a.s. finitely many $p$, and $X$ a r.v. independent of
$\tau$.\\
Then $\mathcal{L}(X,\tau') \ll \mathcal{L}(X, \tau)$.
\end{lm}

This lemma says that changing finitely many values of $\tau$ does
not perturb too much the joint law of $X$ and $\tau$. For instance,
it implies that a process obtained from a Brownian motion by
changing the signs of finitely many excursions has a law absolutely
continuous w.r.t. $\pi$. This is called \textit{'principe de
retournement des excursions'} in [M].\\
\\
\begin{prf}
If $u = (u_1, u_2, ...)$ is an infinite sequence, denote by $u_{p]}$
the finite sequence $(u_1, ..., u_p)$ and by $u_{[p+1}$ the infinite
sequence $(u_{p+1}, u_{p+2}, ...)$. We have $(x,u) = f_p(x, u_{p]},
u_{[p+1})$ for some function $f_p$.\\
We have to show that if $F$ is measurable set such that
$\mathbf{P}[(X,\tau) \in F] = 0$, then $\mathbf{P}[(X,\tau') \in F]
= 0$. So assume $\mathbf{P}[(X,\tau) \in F] = 0$. For $p \in
\mathbb{N}$, since $\tau_{p]}$ takes values in $\{-1, 1\}^p$, we can
write
\begin{equation}
\sum_{\sigma \in \{-1, 1\}^p} \mathbf{P} \left[ f_p (X, \sigma,
\tau_{[p+1}) \in F \ and \ \tau_{p]} = \sigma \right] =
\mathbf{P}\left[(X, \tau) \in F \right] = 0 .
\end{equation}

Using the independence of $\tau_{p]}$ and $(X, \tau_{[p+1})$, this
becomes :
\begin{equation}
\sum _{ \sigma \in \{-1, 1\}^p} 2^{-p} \mathbf{P} \left[ f_p(X,
\sigma, \tau_{[p+1}) \in F \right] = 0 ;
\end{equation}

So for each $p \in \mathbb{N}$ and each $\sigma \in \{-1, 1\}^p$,
the event $\{f_p(X, \sigma, \tau_{[p+1} \in F \}$ is negligible.
Since $\tau'_p(\omega) = \tau_p(\omega)$ for all $p$ larger than
some $N(\omega)$, one has
\begin{eqnarray*}
\mathbf{P} \left[ (X,\tau') \in F \right] &=& \lim_{p\rightarrow
\infty} \mathbf{P} \left[ (X,\tau') \in F \ and \ \tau'_{[p+1} =
\tau_{[p+1} \right]\\
&=& \lim_{p\rightarrow \infty} \sum_{\sigma \in \{-1, 1\}^p}
\mathbf{P} \left[ f_p(X, \sigma, \tau_{[p+1}) \in F \ and \ \tau'_{p]} =
\sigma \ and \ \tau'_{[p+1} = \tau_{[p+1}  \right].
\end{eqnarray*}
This is null because the event ${f_p(X,\sigma,\tau_{[p+1}) \in F}$
is negligible, as shown above.
\end{prf}~
\\
\begin{prop} (mechanism of a Lévy raise)\\
Suppose given the following three r.v. :
\begin{enumerate}
\item \label{l1} $V$, a $\mathbf{W}-$valued r.v., such that $\mathcal{L}(V) \ll
\pi$ ;
\item  \label{l2} $\tau = (\tau_p)_{p \in \mathbb{N}}$, a coin-tossing
independent of $V$ ;
\item  \label{l3} $\tau' = (\tau'_p)_{p \in \mathbb{N}}$, a r.v. valued
in $\{-1, 1\}^\mathbb{N}$, such that the random set ${p \in
\mathbb{N} \qquad \tau'_p (\omega) \neq \tau_p(\omega)}$ is a.s.
finite.
\end{enumerate}

Then there exists a unique $\mathbf{w}-$valued r.v. $U$ such that
\begin{equation}
|U| = V - \mathbf{I} \circ V  \ and \ \mathbf{S}_p \circ U = \tau'_p
\ for \ each \ p.
\end{equation}
It is measurable w.r.t. the $\sigma-$field $\sigma(V, \tau')$ and we
have $\mathcal{L}(U) \ll \pi$ and $\mathbf{T} \circ U = V$. For any $n \geq 0$, we have $U = B^n$ on the event $\{ V = B^{n+1} \ and
\ \tau' = \mathbf{S} \circ B^n \}$.\\
\end{prop}

~\\
\begin{prf}
We start from $\mathcal{L}(V) \ll \pi = \mathcal{L}(B) =
\mathcal{L}(B^1)$. Using Lemma 2 (i) we write $\mathcal{L}(V, V -
\mathbf{I} \circ V) \ll \mathcal{L}(B^1, B^1 - \mathbf{I} \circ
B^1)$. By Lemma 2 (ii), the coin-tossing  $\tau$ (resp. $\mathbf{S}
\circ B$) which is independent of $V$ (resp. $B^1$) can be added on
the left (resp. right), and we obtain $\mathcal{L}(V, V - \mathbf{I}
\circ V, \tau) \ll \mathcal{L}(B^1, B^1 - \mathbf{I} \circ B^1,
\mathbf{S} \circ B)$; by Lemma 1, the right-hand side is
$\mathcal{L}(B^1, |B|, \mathbf{S} \circ B)$. Lemma 3 allows us to
replace $\tau$ by $\tau'$ in the left-hand side, so we finally have
\begin{equation}
\mathcal{L} (V, V - \mathbf{I} \circ V, \tau) \ll \mathcal{L}(B^1,
|B|, \mathbf{S} \circ B)
\end{equation}

Now, we call $\mathbf{W}^+$ the set of non-negative paths and $f :
\mathbf{W}^+ \times \{-1, 1\}^{\mathbb{N}} \longrightarrow
\mathbf{W}$ the measurable function such that $w= f\left(|w|,
\mathbf{S}(w)\right)$. We remark that $B = f\left(|B|, \mathbf{S}
\circ B\right)$ and we define $U= f(V - \mathbf{I} \circ V, \tau')$;
this is the unique r.v. $U$ such that $|U| = V - \mathbf{I} \circ V$
and $\mathbf{S} \circ U = \tau'$. To verify that $\mathcal{L}(U) \ll
\pi$ and $\mathbf{T} \circ U = V$, we apply Lemma 2 (i) to (2) with
the functions $g(x, y, \sigma) = f(y, \sigma)$ and $h(x, y \sigma) =
\left(f(y, \sigma), x\right)$. With $g$ we obtain $\mathcal{L}(U)
\ll \mathcal{L}(B)$, the first claim. With $h$ we obtain
$\mathcal{L}(U, V) \ll \mathcal{L}(B, B^1)$; this implies
$\mathbf{T} \circ U = V$ since the joint law $\mathcal{L}(B), \ (B^1)$
is carried by the graph of $\mathbf{T}$.\\
Last, on the event $\{ V = B^{n+1} and \tau' = \mathbf{s} \circ B^n
\}$, using the definition of $U$ and Lemma 1 we have
\begin{equation}
U=f(V - \mathbf{I} \circ V, \tau') = f(B^{n+1} - \mathbf{I} \circ
B^{n+1}, \mathbf{S} \circ B^n) = f(|B^n|, \mathbf{S} \circ B^n) =
B^n.
\end{equation}
\end{prf} ~
\\
\begin{prop}
Denote by $\mathcal{P}_f(\mathbb{N})$ the set of all finite subsets
of $\mathbb{N}$. Fix $r$ in $\mathbb{N}$. For each $n \leq r$, let be given
$\mathcal{N}^n$, a r.v. with values in $\mathcal{P}_f(\mathbb{N})$,
and $\Sigma^n = (\Sigma_p^n, p \in \mathcal{N}^n)$, a r.v. taking
its values in $\bigcup_{\mathcal{M} \in \mathcal{P}_f(\mathbb{N})}
\{-1, 1\}^{\mathcal{M}}$, such that $\Sigma^n(\omega) \in \{-1,
1\}^{\mathcal{N}^n(\omega)}$.
\\
Starting with $\Gamma^r = B^r$, we can define a sequence $(\Gamma^n)_{n \leq r}$ such that $\Gamma^{n-1}$ is the $\mathbf{W}-$valued r.v. $U$
obtained in Proposition 2 from
\begin{center}
$
V= \Gamma^n, \qquad \tau=\mathbf{S} \circ B'^{n-1}, \qquad \tau'_p =
\left\{%
\begin{array}{ll}
    \Sigma^n_p(\omega) \quad \textrm{ if } \ p \in \mathcal{N}^{n}(\omega)  \\
    \tau_p(\omega) \quad \textrm{ else } \\
\end{array}%
\right.
$
\end{center}
~\\
Then the sequence $(\Gamma^n)_{n \leq r}$ is a $B$-raised Brownian motions of index $r$.
\end{prop}
~\\
\begin{prf}
First, we verify that the $\Gamma^n$ can be constructed stepwise.
Assuming $\Gamma^n$ has already been constructed, has a law absolutely
continuous w.r.t. $\pi$, Proposition 2 applies to $V = \Gamma^{n-1}$ and $\tau = \mathbf{S} \circ
B'^{n-1}$ (they are independent. The r.v. $\Gamma^{n-1} = U$ yielded by Proposition 2 also
satisfies $\mathcal{L}(\Gamma^{n-1}) \ll \pi$, and is measurable in
$\sigma(\mathcal{F}, \tau')$.

The rest of the proof will exhibit a sequence $(\Gamma^n)_{n \geq
r}$ of B-raised motions such that $\Gamma^{J-n} = C^n$. Starting
with $\Gamma^r =B^r$, the other $\Gamma^m$ will be inductively
defined : if $m < r$, suppose $\Gamma^{m+1}$ has been defined, is
$\sigma(B^{m+1})-$measurable, and verifies $\mathcal{L}(\Gamma^{m+1}
\ll \pi$; define $\Gamma^m$ as the r.v. $U$ obtained in Proposition
2 from
\begin{equation}
V= \Gamma^{m+1}, \qquad \tau = \mathbf{S} \circ B^m, \qquad \tau'_p
= \left\{%
\begin{array}{ll}
    \Sigma^{J-m-1}(\omega) \quad \textit{ if } \ p \in \mathcal{N}^{J-m-1}(\omega) \\
    \tau_p(\omega) \quad \textrm{ else } \\
\end{array}%
\right.
\end{equation}
This is possible since $V$ and $\tau$ are independent and
$\mathcal{N}^{J-m-1}$ is a.s. finte; the result $\Gamma^m$ verifies
$\mathcal{L}(\Gamma^n) \ll \pi$ and $\mathbf{T} \circ \Gamma^m =
\Gamma^{m+1}$. To show that $\Gamma^m$ is $\sigma(B^m)-$measurable,
it suffices to show that so is $\tau'$; this may be done separately
on each of the events $\{ J \leq m\}$, $\{ J = m+1 \}$, ..., $\{J =
r \}$, because they form a $\sigma(B^m)-$partition of $\Omega$. On
$\{ J \leq m \}$, we have $\tau' = \mathbf{S} \circ B^m$; this is
$\sigma (B^m)-$measurable. To see what happens for other values of
$j$, introduce $\varphi^n$ and $\psi^n$ such that $\mathcal{N}^n =
\varphi^n (B^{J-n})$ and $\Sigma^n = \psi^n(B^{J-n})$ for $0 \leq n <
k$. For $j \in \{m+1, ..., m+k\}$, we have on $\{J = j\}$
\begin{equation}
\tau' = \psi^n(B^{J-n}) \mathds{1}_{\varphi^n(B^{J-n})} + \tau \left( \mathds{1}_{\Omega} - \mathds{1}_{\psi^n(B^{J-n})} \right)
\end{equation}

This is $\sigma(B^m)-$measurable too. We have established that
$\Gamma^r, ..., \Gamma^0$ exist and form a sequence of $B-$raised
motions; it remains to see that $\Gamma^{J-n} = C^n$.\\
This is done in two steps. Firstly, by induction on $m$, we have
$\Gamma^m = B^m$ on $\{J \leq m\}$ : this holds for $m = r$, and if
it holds for $m+1$, it holds for $m$ too, owing to the last statement
in Proposition 1. Consequently, $\Gamma^m = B^m on \{ J = m\}$, that
is $\Gamma^J = B^J = C^0$. Secondly, to proceed by induction on $n$,
we will assume that $\Gamma^{J-n} = C^n$ for some $n \geq 0$, and
show $\Gamma^{J-n-1} = C^{n+1}$. It suffices to show this equality
on the event $\{J = j \}$; on this event, using the definition of
$\Gamma^m$ with $m = j-n-1$ and the inequality $m = j-n-1 < j$, the
r.v. $\Gamma^{J-n-1}$ satisfies both $\mathbf{T}(\Gamma^{j-n-1}) = \Gamma^{j-n} = C^n =
\mathbf{T}(C^{n+1}$) and
\begin{equation}
\mathbf{S} \circ \Gamma^{j-n-1} = \left\{%
\begin{array}{ll}
    \Sigma^{n+1} \quad \textrm{ on } \ \mathcal{N}^n \\
    \mathbf{S} \circ B^{j-n-1}  \quad \textrm{ else } \\
\end{array}%
\right.
= \mathbf{S} \circ C^{n+1}.
\end{equation}
These two equalities entail $\Gamma^{J-n-1} = C^{n+1}$ a.s. on $\{J
= j\}$.

\end{prf}~
\\
\begin{lm}
Let $(j,k) \in \mathbb{N}^2$ and $Q$ and $R$ be two r.v. such that
$k \leq j-1$ and  $0\leq Q \leq R$. On the event $\left\{\forall n
\in \{k, ..., j-1\} Z \circ B^n \cap (Q, R) = \emptyset \right\}$
that the first iterates of $B$ do not vanish between $Q$ and $R$,
there exists a (random) isometry $i : \mathbb{R} \longrightarrow
\mathbb{R}$ such that $B^j = i \circ B^k$ on the interval $(Q, R)$.
\end{lm}
~\\
\begin{prf}
By induction, it suffices to show that if $B^{j-1}$ does not vanish
on the interval $(Q,R)$, then $B^j = i \circ B^{j-1}$ on $(Q,R)$,
for some random isometry $i$. This is just Lemma 5 with $j = k + 1$
and $B^{j-1}$
instead of $B^k$, so we may suppose that $j=1$.\\
On the event $\{Q = R\}$, the result is trivial. On $\{Q < R\} \cap
\{Z \circ B \cap (Q,R) = \emptyset \}$, the local time $L$ is
constant on $[Q,R]$ because its support is $Z \circ B$, and the sign
of $B$ is constant on $(Q,R)$; so $B^1 = |B| - L = i(B)$ on $(Q,R)$,
where $i$ is the random isometry $x \mapsto x sgn
\left(B_{(Q+R)/2}\right) - L_{(Q+R)/2}$.
\end{prf}~
\\
\begin{nota}
For $w \in \mathbf{W}$, the p-th excursion interval $e_p(w)$ was
defined earlier; the number $\mathbf{h}_p(w) = \max\limits_{s \in
e_p(w)} |w(s)|$ will be called the height of the corresponding
excursion.
\end{nota}

\begin{lm}
Let $X$ be a process whose law is absolutely continuous w.r.t.
Wiener measure. Almost surely,
\begin{itemize}
\item $\lim\limits_{p\rightarrow\infty} \mathbf{h}_p(X)\mathds{1}_{\{e_p(X)\subset[0,t]\}} = 0$;
\item $\sum\limits_{p \in \mathbb{N}} \mathbf{h}_p(X) \mathds{1}_{\{e_p(X)\subset[0,t]\}} =
\infty$;
\item the set $\left\{ \sum\limits_{p \in \mathcal{M}} \mathbf{h}_p(X) \mathds{1}_{\{e_p(X)\subset[0,t]\}}, \mathcal{M} \in \mathcal{P}_f(\mathbb{N})\right\} \ is \ dense \ in \ [0,
\infty)$;
\item between any two different excursions of $X$, there exists a
third one, with height smaller than any given random variable $\eta
> 0$.
\end{itemize}
\end{lm}

~\\
\begin{prf}
By a change of probability, we may suppose that $X$ is a Brownian
motion. It is known (see Exercise (VI.1.19) of [RY]) that when $\eta
\rightarrow 0^+$, the number $\sum\limits_p
\mathds{1}_{\{e_p(X)\subset[0,t]\}} \mathds{1}_{\{\mathbf{h}_p(X) >
\eta\}}$ of downcrossings of the interval $[0, \eta]$ by $|X|$
before $t$ is a.s. equivalent to $\eta^{-1}L_t$, where $L_t$ is the
local time of $X$
at $0$. This easily implies (i) and (ii), wherefrom (iii) follows.\\
Last, between any two excursions of $X$ there are infinitely many
other ones (because $X$ has no isolated zeroes) and, by (i), only
finitely many with heights above $\eta$, whence (iv).
\end{prf}~
\\
\begin{nota}
An excursion whose interval is included in $[0,t]$ will be called a
\textit{$t-$excursion}.
\end{nota}

 It remains to describe the $\mathcal{N}^n$ and $\Sigma^n$, i.e., to choose the
signs of finitely many excursions when Lévy-raising from $\Gamma^n$ to
$\Gamma^{n-1}$. This will be done soon; we first need some notation and a
lemma.

\begin{nota}
If $e'$ and $e''$ are two excursions of a path (or of a process),
$e' \prec e''$ means that $e'$ is anterior to $e''$ : $s' < s''$ for
all $s' \in e'$ and $s'' \in e''$.\\
For an excursion $e$ of $w$, we denote by $i_w e:= \inf \{ w_s; s\in [0, d_e] \}$.
\end{nota}

\begin{df}
An excursion $e$ of a path $w \in \mathbf{W}$ is said to be
\textit{tall} if it is positive (this implies that the process
$\mathbf{I}w$ remains constant during $e$); and if for any excursion $e'$
of $w$ such that  $i_we' =
i_we$ and higher than $e$, then $e'=e$.
\\
Formally, $e$ is tall if it is positive and if
\begin{equation}
\max\left( w(s); s \geq 0, (\mathbf{I}w)(s) = i_we \right) = \max
\left(w(s); s \in e \right).
\end{equation}
\end{df}~
\\
\begin{lm}
Let $\eta$ be a positive number, $m \geq 1$ be an integer and $w \in
\mathbf{W}$ a path. Let $e_1, \dots, e_{m+1}$ be $m+1$ different
$t-$excursions of $w$, numbered in chronological order : $e_1 \prec
\dots \prec e_{m+1}$; call $h_1, \dots, h_{m+1}$ their respective
heights. Let $f_1, \dots, f_p$ denote all excursions of $w$ which
are anterior to $e_{m+1}$ and whose heights are $\geq \min( \eta,
h_1, \dots, h_{m+1})$, numbered in reverse chronological order : let $g_1, \dots, g_p$ be $p$ excursions
of $w$ verifying $f_p \prec g_p
\prec \dots \prec f_1 \prec g_1 \prec e_{m+1}$.\\
Suppose that
\begin{itemize}
\item the excursion $e_{m+1}$ is negative, and all $t-$excursions higher
than $e_{m+1}$ are positive;
\item the excursions $f_1, \dots, f_p$ are positive;
\item the excursions $g_1, \dots, g_p$ are negative; and every
negative excursion anterior to $g_p$ is smaller than $g_q$.
\end{itemize}
We call the $g_i$'s the plug-excursions, and the $e_j$'s the protected excursions.
Then $e_1, \dots, e_m$ are tall, and $|i_we_1| < |i_we_2| < \dots <
|i_we_m| < \eta$.
\end{lm}

~\\
\begin{prf}
Firstly, $|i_wf_1| < \eta$ because $f_1 \prec e_{m+1}$ and any
excursion anterior to $e_{m+1}$ and having height $\geq \eta$ is one
of the $f_q$,
hence positive.\\
Secondly, for $1 \leq q \leq p$, the excursion $g_q$ is negative and
higher than any negative excursion, anterior to it; so $\mathbf{I}w$ is not
constant during $g_q$, and consequently we have
\begin{equation}
|i_wf_p| < |i_wf_{p-1}| < \dots < |i_wf_1| < \ height \ of \ g_1,
\end{equation}
where each $<$ sign is due to $\mathbf{I}w$ varying on the
corresponding
$g_q$.\\
thirdly, combining (20) with $|i_wf_1| < \eta$ (first step), and
noticing that, by definition of the $f_q$, $(e_1, \dots, e_m)$ is a
sub-sequence of $(f_p, \dots, f_1)$, we obtain
\begin{equation}
|i_we_1| < \dots < |i_we_m| < \eta.
\end{equation}
Last, it remains to establish that $e_l$ is tall for $1 \leq l \leq
m$. Let $e'$ denote a positive excursion of $w$ with height $h' \geq
h_l$ and such that $i_we' = i_we_l$. From (13), we have $|i_we'| =
|i_we_l| < height \ of \ g_1$; so $e'$ is anterior to $g_1$ and a
fortiori anterior to $e_{m+1}$. As $h' \geq h_l$, $e'$ must be one
of the $f_q$ (see their definition). But $e_l$ is also one of the
$f_q$ and, due to (13), all $i_wf_q$ are different; so $e'=e_l$. This
means that $e_l$ is tall.
\end{prf}~
\\
In the proof of Lemma 6, the negative excursions $g_q$ are used to
separate the $f_q$ from each other. Yet, in the end, we are not
interested in the behavior of all $f_q$ but only in the $e_l$. It is
possible to replace this lemma with a variant, where $2m$ excursions
(instead of $p$ ones, the $g_q$) are made negative, each $e_l$ being
flanked by two of them.\\
\\
\begin{lm}
Let $X$ be a process with law absolutely continuous w.r.t. $\pi$,
and $E$ a tall excursion of $\mathbf{T} \circ X$ with height $H$.
There exists an excursion of $X$, with interval $\{s; (\mathbf{I}
\circ \mathbf{T} \circ X)(s) = i_{\mathbf{T} \circ X} E \}$, and
with height $H + |i_{\mathbf{T} \circ X} E|$.
\end{lm}

~\\
\begin{prf}
First, recall a.s., Brownian motion $B$ does not reach its current
minimum $\mathbf{I} \circ B$ in the interior of a time-interval
where $\mathbf{I} \circ B$ is constant. (This is a consequence of
($\mathbf{I} \circ B)(s) < 0$ for $s > 0$ and of the Markov property
at
the first time that $B =\mathbf{I} \circ B$ after some rational).\\
Put $Y = \mathbf{T} \circ X$ and call $F$ the interval $\{s \geq 0;
(\mathbf{I} \circ Y)(s) = i_XE\}$; $Y$ reaches its current minimum
$\mathbf{I} \circ Y$ at both endpoints of $F$ but not in the
interior of $F$ (see above). Since $|X| = Y - \mathbf{I} \circ Y$ by
Lemma 1, we have that $F$ is the support of some excursion of $X$.
The height of that excursion is
\begin{eqnarray*}
\max (|X_s| ; s \geq 0, &and& ( \mathbf{I} \circ Y)(s) = i_Y E) \\
&=& \max \left(Y_s - (\mathbf{I} \circ Y)(s); s \geq 0, (\mathbf{I}
\circ Y)(s) = i_Y E
\right)\\
&=& \max \left(Y_s ; s \geq 0, (\mathbf{I} \circ Y)(s) = i_Y E
\right) - i_Y E\\
&=& \max \left(Y_s ; s \in E \right) - i_Y E \qquad \ because \ E \ is \ tall \\
&=& H + |i_Y E|.
\end{eqnarray*}
\end{prf}~
\\
\begin{lm}
Let $(0,\overrightarrow{i},\overrightarrow{j})$ be an orthonormal
basis of the plan in which we represent paths. Let $\tau^{a+}_b$ be
the vertical translation of vector $(b-a)\overrightarrow{j}$ and
$\tau^{a-}_b$ the reflection along the horizontal axis of equation :
$\dis y=\frac{a+b}{2}$.\\
Consider $(t,k,p)\in \mathbb{R}^+_*\times\mathbb{N}^2$ such that
$w^k_t=a$ and $w^{k+p}_t=b$ and denote $\gamma_t$ the first time
posterior to $t$ when at least one of the iterated Lévy transforms
$w^s$, $k\leq s\leq k+p-1$, vanishes. Then we have :\\
\begin{center}
$w^{k+p}_{\emph{\textbf{|}}[t,\gamma_t]}=
\left\{%
\begin{array}{ll}
    \tau^{a+}_b\;o\,w^k_{\emph{\textbf{|}}[t,\gamma_t]}\quad\textrm{ if }\qquad\underset{i=k}{\overset{k+p-1}{\Pi}}w^i_t>0 \\
    \tau^{a-}_b\;o\,w^{-k}_{\emph{\textbf{|}}[t,\gamma_t]}\quad\textrm{ else } \\
\end{array}%
\right.$ \end{center}~\\
We will denote $\tau^k_{k+p}(w)$ the plan transformation, which
transforms $w^k_{\emph{\textbf{|}}[t,\gamma_{t}]}$ in
$w^{k+p}_{\emph{\textbf{|}}[t,\gamma_{t}]}$.
\end{lm}
~\\
\begin{prf}
It is an immediate consequence of Tanaka's Lemma, when $p=1$. \\
In general case, we break up the displacement $\tau$ which
transforms $w^k_{\emph{\textbf{|}}[t,\gamma_{t}]}$ in
$w^{k+p}_{\emph{\textbf{|}}[t,\gamma_{t}]}$ under the form
$\tau=\tau_p\circ\tau_{p-1}\circ...\circ\tau_1$ where $\tau_i$
transforms $w^{k+i-1}_{\emph{\textbf{|}}[t,\gamma_{t}]}$ in
$w^{k+i}_{\emph{\textbf{|}}[t,\gamma_{t}]}$. From the preceding
remark, each $\tau_i$ is a vertical translation or a reflection
along an horizontal axis, according to the sign of $w^{k+i-1}_t$.
Then we deduce the claim.
\end{prf}
~\\
To construct the desired process $\Gamma$, we will proceed by induction on discretized time, and so we will perform, from the level $r$, two types of raisings.\\
In a first type, the so-called horizontal raisings, at a level when the path vanishes on$ [t_d,t_{d+1}]$, we protect the material furnished by the induction hypothesis on $[0, t_d]$, namely $\widetilde{\Gamma}$. And we prepare the path on $[t_d, t_{d+1}]$ to give it the form it ought to have at this level for being near to $\varphi$ on the interval when $\Gamma$ is near to $\varphi$ on $[0,t_d]$. \\
So we make positive the significative excursions of the path called here the protected excursions, and insert between them small excursions called the plug-excursions (see lemma 6). And on $[t_d, t_{d+1}]$, we prepare excursions , the building ones, which we protect, and they will act, one by one, during a succession of horizontal raisings, to give the path the previewed form (see Lemma 13 and Proposition 14), at the condition the path, after that, will not vanished on $[t_d, t_{d+1}]$.\\
In a second type, the so called vertical raisings, we give anew to the protected excursions the signs they have before the horizontal raisings. Then we get up while the path doesn't vanish on $[t_d, t_{d+1}]$. In fact we must distinguish the last raising of a succession of horizontal raisings, the so called terminal horizontal raising, when we leave to protect the protected excursions and give them the good signs.\\
It is important to know the real level of the path, ie. the level without the horizontal raisings. Precisely, we define the r.v. inductively :
\begin{eqnarray*}
RL_r(w) = r \\
and \ for \ all \ integer \ n \leq r \\
RL_{n-1}(w) =
\begin{cases}
RL_n(w) \quad \textit{if the step } n-1 \rightarrow n \textit{ corresponds to an horizontal raising} \\
RL_n(w) - 1 \quad \textit{if it corresponds to a vertical raising.}
\end{cases}
\end{eqnarray*}

For our needs, we will call map-excursion, or simply excursion, each map $e : \mathbb{R}^+
\to \mathbb{R}$ whose support is a not empty segment and which
doesn't vanish at any point of the interior of the support. In
particular, for $w\in W$ and $t>0$, we will call excursion
straddling $t$, and denote it
by : $e_t(w)$, the map so defined :\\
\\
$e_t(w):\mathbb{R}^+\to\mathbb{R},\,\forall
u\in\mathbb{R}^+,\,e_t(w)(u)=
\left\{%
\begin{array}{cl}
    0&\textrm{ if }u\in[0,g_t(w)]\cup[d_t(w),+\infty[ \\
    w_u&\textrm{ else } \\
\end{array}%
\right.$\\\\
We will introduce the map $de_t(w) :
\mathbb{R}^+\to\mathbb{R}$ defined by\\\\
$\forall u\in\mathbb{R}^+,\;de_t(w)(u)=
\left\{%
\begin{array}{cl}
    0 &\textrm{if }u\in[0,t]\cup[R_t(w),+\infty[\\
    &\textrm{where }R_t(w)=\sup\{u>t,\forall s\in]t,u],w_s\neq w_t\} \\
    w_u-w_t &\textrm{else } \\
\end{array}%
\right.$\\\\
and we call it a differential excursion of $w$ whenever its support is not empty, ie.
\begin{equation*}
\exists \varepsilon >0, \ \forall u \in (t, t+\varepsilon), \ w_u \neq w_{\varepsilon}
\end{equation*}
We denote them $e_t$ and $de_t$ when there is no ambiguity, and
$e^s_t$ and $de^s_t$ in the case of excursions of $w^s$.
\\

\begin{lm}
Let $w\in W$ and $e$ a negative $m-$excursion of $\mathbf{T}w$,
lower than all preceding it. Let $\gamma$ be its beginning and
$\delta$ its end. We
set :\\
\begin{equation*}
    \gamma_1=\underset{[0,\gamma]}{\arg\min\,\mathbf{T}w},\,\gamma_2=\inf\{t\in
    supp(e);e(t)= \mathbf{T}w_{\gamma_1}\},\,\gamma_3=\underset{[\gamma,\delta]}{\arg\min\,e}
\end{equation*}
Then :\\
\[
\left\{
\begin{array}{l}
    de_{\gamma_1}\,\textrm{ coincides with an excursion of }|w|,\textrm{ and its support is }[\gamma_1,\gamma_2]\\
    de_{\gamma_3}\,\textrm{ coincides with an excursion of }|w|,\textrm{ which begins at }\gamma_3\textrm{ and whose }\\
    \textrm{ support contains }[\gamma_3,\delta]\\
\end{array}
\right.
\]
\\
Furthermore, $\forall u\in[\gamma_2,\gamma_3]$, $de_u$ coincides
with an excursion of $|w|$, if, and only if :
\[\left\{
\begin{array}{c}
    de_u\,\textrm{ is a positive excursion } \\
    e(u)=\inf\{e(t),\,t\in[\gamma,u]\} \\
\end{array}
\right.
\]
It is the case in particular when $de_u$ is the first positive
excursion of the form $de_v$, $v\in [\gamma_2,\gamma_3]$ to overflow
a given value.
\end{lm}
~\\
\begin{prf}\\
From Tanaka's formula :\\
\begin{equation*}
|w_t| = \mathbf{T}w_t+\sup\{-\mathbf{T}w_u, u\in [0,t]\}
\end{equation*}
Therefore,
\begin{equation*}
|w_{\gamma_1}|=\mathbf{T}w_{\gamma_1}-\mathbf{T}w_{\gamma_1}=0,
\end{equation*}
while, for all $t>\gamma_1$, sufficiently small :
\begin{equation*}
\mathbf{T}w_t > \mathbf{T}w_{\gamma_1}.
\end{equation*}
So, $de_{\gamma_1}$ is a positive excursion of
$|w|$ which ends at $\gamma_2$.\\
In the same way, $\mathbf{T}w_t > \mathbf{T}w_{\gamma_3}$, for all
$t\in [\gamma_3,\delta]$, therefore $de_{\gamma_3}$ is an excursion
of $|w|$ beginning at $\gamma_3$ whose support
contains $[\gamma_3,\delta]$.\\
Let $e'$ be an excursion of $w$ with support included in
$[\gamma_2,\gamma_3]$. Its beginning $u$, and its end $v$ verify :\\
\begin{equation*}
u=\underset{[0,v[}\arg\min\, \mathbf{T}w\qquad\textrm{and}\qquad
v=\underset{[0,v]}\arg\min\, \mathbf{T}w.
\end{equation*}

\noindent So we deduce : \qquad $de_u=|e'|$.\\
Reciprocally, let $u\in [\gamma_2,\gamma_3]$ such that $de_u$ is a
positive excursion and\\
$u=\underset{[0,u]}\arg\min\, \mathbf{T}w$.
\\
Then, $u=\underset{[0,v[}\arg\min\, \mathbf{T}w$, where $v$ is the
end of $de_u$, because $de_u$ is positive.
\\
Thus, $w_u = w_v = 0$, and for all $t\in ]u,v[$,
$w_t \neq 0$.\\
Consequently, $de_u$ is an excursion of $|w|$.\\
Let $h>0$ be such that there exists $u\in [\gamma_2,\gamma_3]$
verifying $de_u$ is the first positive excursion of the form $de_v$,
$v\in [\gamma_2,\gamma_3]$, whose height overflows $h$. Then, for
all $v<u$, the support of $de_v$ can't contain this of $de_u$
without denying the minimality of $u$.
\end{prf}~
\\
As an immediate consequence, we observe : \\
\begin{cor}
The excursions of $w$ coincide with the positive differential excursions of $\mathbf{T}w$, beginning at $\arg \min_{[0,t]}(\mathbf{T}w)$ for all $t \in [0, +\infty[$.
\end{cor}

\newpage
\section{Density of orbits}
\noindent In this paragraph, we want the raised path to approach the
map $\varphi$ uniformly on $[0,T]$. Precisely :\\
\textit{Whatever $\varepsilon$ strictly positive, and $\varphi \in
W_{|_{[0,T]}}$, there exists $\Gamma$ a $B-$raised
Brownian motion such that :
\begin{equation*}
    \mathbf{P}\left(\|\Gamma_{|_{[0,T]}}-\varphi_{|_{[0,T]}}\|_\infty<\varepsilon\right)>1-\varepsilon
\end{equation*}}
\\\\
We consider a modulus of uniform continuity $\alpha_0$ associated to
$\dis\left(\frac{\varepsilon}{4},\varphi,[0,T]\right)$ and a real
number $\alpha_1$ such that $\dis
P(A_{0\,\varepsilon})>1-\frac{\varepsilon}{2}$ where
\begin{equation*}
    A_{0\,\varepsilon}=\left[\sup\{|B_t-B_u|,\,(t,u)\in[0,T]^2\textrm{
    and
    }|t-u|<\alpha_1\}<\frac{\varepsilon}{2}\right]\,,
\end{equation*}
then we set $\alpha:=\min(\alpha_0,\alpha_1)$,
$d_0:=\left[\frac{T}{\alpha}\right]+1$, and for all $d\in
\mathbb{N}$,
$t_d=(d\alpha)\wedge T$.\\
We set again, for all integer $d\in [1,...,d_0]$,
\begin{equation*}
A^d_\varepsilon:=\left[\sup\{|B_t-B_u|,\,(t,u)\in[t_d,T]^2,\,|t-u|<\alpha\}<\frac{\varepsilon}{2}\right]\,.
\end{equation*}
Our aim is to show, by induction on $d$, the following property
$\mathcal{P}_d$ : " For all $\varepsilon > 0$, there exists an
integer $r_d$ and $\Gamma$ a $B-$raised Brownian
motion of index $r_d$ such that :
\begin{equation*}
    \mathbf{P}\left(\left[\|\Gamma_{|_{[0,t_d]}}-\varphi_{|_{[0,t_d]}}\|_\infty<\varepsilon\right]\cap\left[|\Gamma_{t_d}-\varphi(t_d)|<\varepsilon_1\right]\cap
    A_\varepsilon^d\right)>1-\varepsilon\left(1+\frac{d}{d_0}\right)\;"
\end{equation*}
Notice that $\mathcal{P}_0$ immediately yields from the choice of
$\varphi$ which vanishes at $0$. We suppose now $\mathcal{P}_d$ true. We are going to
apply this hypothesis to the Brownian motion $B^{s_0}$, for an
integer $s_0$ which, as the real number
$\varepsilon_{1}$, will be later specified.\\
As
$A^d_\varepsilon\subset\left[\sup\{|B_t-B_u|,\,(t,u)\in[t_d,t_{d+1}]^2\}<\frac{\varepsilon}{2}\right]\cap
A^{d+1}_\varepsilon$, and from the independence of the increments of
Brownian motion, we can deduce the existence of a disjointed sum
$\widetilde{\Gamma}$ of $B^{s_0}-$raised Brownian motions of index
$r_d$ such that :
\begin{equation*}
    \mathbf{P}(A_0^\varepsilon)>1-\varepsilon\left(1+\frac{d}{d_0}\right)\;,
\end{equation*}
where
\begin{eqnarray*}
    A_0^\varepsilon&:=&\left[\|\widetilde{\Gamma}_{|_{[0,t_d]}}-\varphi_{|_{[0,t_d]}}\|_\infty<\frac{\varepsilon}{2}\right]\cap\left[|\widetilde{\Gamma}_{t_d}-\varphi(t_d)|<\varepsilon_1\right]\quad...\\
    &&...\quad\cap\left[\sup\{|B_t^{s_0+r_d}-B_u^{s_0+r_d}|,\,(t,u)\in[t_d,t_{d+1}]^2\}<\frac{\varepsilon}{2}\right]\cap A^{d+1}_\varepsilon
\end{eqnarray*}
(It suffices to apply $\mathcal{P}_d$ with $\frac{\varepsilon}{2}$ instead of $\varepsilon$).\\

We will denote : $\forall
i\in\mathbb{N},\,\widetilde{\Gamma}^i=\widetilde{w}^i$. By
definition, $\forall i>r_d,\,\tilde{w}^i=w^{s_0+i}$.
\\
From the theorem of density of zeroes ([M]), there exists a.s. an
integer $\ell$ such that $\widetilde{w}^\ell$ vanishes at
least one time on $[t_d, t_{d+1}]$.\\
Let $L(w)$ be the smallest of these integers $\ell$. $L$ is a r.v.
almost surely finite. Then there exists an integer $\ell_0$ which we
will choose $>r_d$ such that :
\begin{equation*}
    \mathbf{P}(A_1^\varepsilon)>1-\varepsilon\left(1+\frac{d}{d_0}\right)-\frac{\varepsilon}{2d_0}\quad\textrm{où}\quad
    A_1^\varepsilon:=A_0^\varepsilon\cap[L\leq\ell_0]\,.
\end{equation*}
We set : $r=s_0 + \max (r_d, l_0)$.

\begin{df}
The protecting tree of $\widetilde{\Gamma}$.\\
Let $\mathcal{N}$ be a set of excursions of $w_{|[0, d_t(w)]}$, $w\in \mathbf{W}$, and $\varepsilon >0$. We denote by $\mathbf{T}_w^{t, \varepsilon}(\mathcal{N})$ the set of excursions of $\mathbf{T}w_{|[0, d_t(\mathbf{T}w)]}$ so defined : \\
$e' \in \mathbf{T}_w^{t, \varepsilon}(\mathcal{N})$ iff $e'$ is an excursion of $\mathbf{T}w_{|[0, d_t(\mathbf{T}w)]}$ which appears in the differential excursions $dE$ of $\mathbf{T}w$ corresponding to an excursion $E \in \mathcal{N}$, and $h(e') > \frac{\varepsilon}{4}$ (we call them excursions of the first type) or the support of $e'$ contains $\arg\max$ of the differential excursion $dE$, we call them excursions of the second type, and all the excursions $e''$ of $\mathbf{T}w$ of height belonging to $[h(e'), h(e^{'}_{-})]$, $e^{'}_-$ being the preceding excursion of $\mathbf{T}w$ of height $> \frac{\varepsilon}{4}$, with $e^{'}_- < e''<e'$, we call these excursions $e''$ the third type excursions.
\end{df}~
\\
\begin{lm}
When $\mathcal{N}$ is a finite set, so is $\mathbf{T}_w^{t, \varepsilon}(\mathcal{N})$.
\end{lm}
\begin{prf}
We remark firstly that the number of distinct differential excursions of $\mathbf{T}w$ corresponding to the elements of $\mathcal{N}$ is finite, equal to the cardinal of $\mathcal{N}$. The number of excursions of the first type is finite because their height is greater than $\frac{\varepsilon}{4}$. The number of excursions of the second type is finite following the first remark. And for each excursion of the second type, the number of excursions of the third type is finite too.
\end{prf}

Then we call protecting tree of $\widetilde{\Gamma}$, the tree constituted by :
\begin{itemize}
\item at the $0$-generation : the elements of $\widetilde{\mathcal{N}}^0$, the set of $d_t(\widetilde{w}^0)$-excursions of height $> \varepsilon$.
\item and, for all $n \in \mathbb{N}$, if we denote by $\widetilde{\mathcal{N}}^n$, the set of protected excursions of $\widetilde{w}^n$ whose elements, given in chronological order, constitute the $n^{th}$ generation of the tree, then at the $(n+1)^{th}$ generation, the elements of $\widetilde{\mathcal{N}}^{n+1} := \mathbf{T}_{\widetilde{w}^n}^{t_d, \varepsilon / 4^n}\left(\widetilde{\mathcal{N}}^n\right)$.
    As we have ordered each $\widetilde{\mathcal{N}}^n$, which is now a finite sequence of excursions of $\widetilde{w}^n$, we can consider $\widetilde{\Sigma}^n$ the finite sequence of their signs.
    And our next task is to definite $\mathcal{N}'^n$ the sequence of protected $t_d$-excursions of $\Gamma^n$ and $\Sigma '^n$ the associated sequence of their signs.
    We have defined the $\widetilde{\mathcal{N}}^n$'s by getting down in the iterations. We will define the $\mathcal{N}'^n$'s by getting up from level $r$.
    We put $\mathcal{N}'^r(\omega) = \widetilde{\mathcal{N}}^{r-s_0}(\omega)$, and  $\Sigma '^r(\omega) = \widetilde{\Sigma}^{r-s_0}(\omega)$.\\
    If $\Gamma^{n+1} \longrightarrow \Gamma^n$ is an horizontal raising, then $\mathcal{N}^{0,n+1}$ is constituted by the family of the excursions of $\Gamma^n$ whose absolute values have the same $\arg \max$ as the excursions of $\mathcal{N}^{0,n+2}$. Then we proceed as in lemma 6, considering the elements of $\mathcal{N}^{0,n+1}$ as the $e_i$, $1 \leq i \leq m$.
    We have now to define the r.v. $\eta$ :
    let $\eta_i$ be the minimal distance between the heights of distinct protected excursions of $\widetilde{w}^i$ for $i \in \{0, \dots, r-s_0\}$ and $\eta_m = \bigwedge \limits_{i=0}^{r-s_0} \eta_i, \widetilde{\varepsilon}_i = \varepsilon \wedge \eta_m \frac{1}{4^{i+1}}$, for $i \in \{0, \dots, r-s_0\}$, and $\varepsilon_i = \widetilde{\varepsilon}_{RL(i) - s_0}\frac{1}{4^{i'-i}}$, where $i' = \sup \{j>i/RL(j) = RL(i)\}$.
    Then we choose for $\eta$ the r.v. $\varepsilon_{RL(n) - s_0}$.
    And we (can) define the set $\mathcal{N}'^n$ as the set of elements the $f_j$'s and the $g_j$'s. And $\Sigma'^n$ give sign $-1$ to the $g_j$'s and $+1$ to the $f_j$'s if $\Gamma^{n+1} \longrightarrow \Gamma^n$ is not a terminal horizontal raising.
\end{itemize}
    In the case of a terminal horizontal raising, we put simply $\mathcal{N}'^{n+1} = \mathcal{N}^{0,n+1}$.\\
    If $\Gamma^{n+1} \longrightarrow \Gamma^{n}$ is a vertical raising, $\mathcal{N}^{n+1}$ is constituted by the family of the excursions of $\Gamma^n$ corresponding to the positive differential excursions of $\Gamma^{n+1}$ whose support encounters at least the support of an element of $\mathcal{N}^{0, n+2}$ and beginning at an $\arg \min$ of $\Gamma^{n+1}$. In these two cases, to specify $\Sigma'^{n+1}$, we need the following lemma. Let us call argext of an excursion the $\arg\max$ (resp. $\arg\min$) of this excursion if it is positive (resp. negative).
\begin{lm}
If for all $i$ from level $r$ to $n$ (with $r\geq n$) the elements of $\mathcal{N}^{0,i}$ and $\widetilde{\mathcal{N}}^{RL(i) - s_0}$ have the same argext, which implies $\left|\mathcal{N}^{0,i}\right| = \left|\widetilde{\mathcal{N}}^{RL(i) - s_0}\right|$, and the same signs, then the elements of $\mathcal{N}^{0,n-1}$ and $\widetilde{\mathcal{N}}^{RL(n-1) - s_0}$ have the same argext, hence these sets have the same cardinality, and the same hierarchy, ie. the order of the heights between the respective excursions of each set is the same. So we put $\Sigma^{0,n-1} = \widetilde{\Sigma}^{RL(n-1)-s_0}$.
\end{lm}

\begin{prf}
The distinction between $\Gamma$ and $\widetilde{\Gamma}$ is due to the introduction of horizontal raises.
If $\Gamma^{n+1} \longrightarrow \Gamma^n$ is an horizontal raising, each protected excursion gains in height the height of a plug excursion, and its support enlarges.
So, during such a raising, the "error" between $\Gamma$ and $\widetilde{\Gamma}$ increased of the height of a plug excursion. For the time, in this lemma, the "error" means the maximal distance between the heights of corresponding excursions.
If, on the other side, $\Gamma^{n+1} \longrightarrow \Gamma^n$ is a vertical raising, a protected excursion of $\Gamma^n$ is constituted by at most two protected excursions of $\Gamma^{n+1}$. So during such a raise, the error between $\Gamma$ and $\widetilde{\Gamma}$ doubles.
And it is easy, by means of our choice of the heights of the plug excursions, to show that this error is always majorized by $\varepsilon \wedge \eta$. Hence the lemma follows immediately.
\end{prf}
Now we can consider (see Proposition 2) that $\Gamma_{|[0,t_d]}$ is correctly defined : for once $RL_w(n) = s_0$, we end the raises, if necessary, by horizontal ones : protecting all the material acquainted.
But is $\Gamma$ near from $\widetilde{\Gamma}$ ? To answer this question, we have to consider the error between $\Gamma$ and $\widetilde{\Gamma}$ now, as being $\| \Gamma_{|[0,t_d]} - \widetilde{\Gamma}_{|[0,t_d]} \|_{\infty}$.

We have just built $\Gamma_{|[0,t_d]}$. We have now to control : $\| \Gamma^N_{|[0,t_d]} - \widetilde{\Gamma}^0_{|[0,t_d]} \|_{\infty}$, where $N = \sup \left\{n\in \{0,\dots,r\} / RL(n) - s_0 = 0 \right\}$ (with $N= - \infty$ if the set is empty). $N$ is a r.v.

\begin{lm}
$\| \Gamma^N_{|[0,t_d]} - \widetilde{\Gamma}^0_{|[0,t_d]} \|_{\infty} \leq 2\varepsilon$, on the event $\left[N\geq 0\right]$.
\end{lm}

\begin{prf}
Let us first remark, from the preceding lemma, that, on the event $\left[N \geq 0 \right]$, we have protected excursions with the same argext, the same signs, and heights near at $\varepsilon$. The supports of the protected excursions of $\Gamma^N$ containing the supports of the corresponding protected excursions of $\widetilde{\Gamma}^0$.
Furthermore, on the difference of their supports, $\Gamma^N$ and $\widetilde{\Gamma}^0$ differ from at most $2\varepsilon$, and likely outside the union of their supports.

\end{prf}

The purpose of the following lemma is to prepare, at level $s$, when
the iterated Brownian motion vanishes on $]t_d,t_{d+1}[$, the
excursions which will allow the correctly raised path to approach
$\varphi$ at level $0$ on $]t_d,t_{d+1}[$.
\\

\begin{lm} Full planing.
\\
Let $w$ belong to $W$, and $t, \ t', \ \varepsilon' \ \in
\mathbb{R}^+_*$ be such that $t < t'$. We suppose there is no
interval in which $w$ is constant, and $w$ vanishes in $(t, t')$.\\
The following r.v. are functionals of $|w|$ :
\begin{itemize}
    \item [$\bullet$] $t_0(w) = g_{\inf \{ s > d_t; |w_s| \geq
    \varepsilon'\}} \wedge t^{\prime}$, with $\inf \emptyset = + \infty$.
    \item [$\bullet$] $\forall n \in \mathbb{N}$, while
$t_n < t'$, we set :\\
$t_{n+1}:=
\left\{%
\begin{array}{ll}
    \inf\{u\in[t_n;\arg\max|de_{t_n}|[, h(de_u) > \varepsilon'
    \textrm{ and sgn}(de_u) = -\textrm{sgn}(de_{t_n})\}\\
    \textrm{if this set is not empty,} \\
    \textrm{else :}\\
    \arg\max|de_{t_n}| \wedge t'
\end{array}%
\right.
$\\
\end{itemize}
The sequence $(t_n)_{n\in \mathbb{N}}$ is strictly increasing and
finite. Let $1+K(w)$ be its cardinality.
\end{lm}
~\\
\begin{prf}
By construction, the sequence $(t_n)$ is strictly increasing and
lower than $t'$. Suppose the number of its terms is infinite. In this
case, it would admit a limit $t_\ast\leq t'$, and the oscillation of
$w$ at $t_\ast$ would be infinite, so contradicting the continuity
of $w$. Then $(t_n)_n$
is finite.\\
The measurability and the finiteness of $K$ are immediate.
\end{prf}~
\\
Let us remark that this Lemma gives us the possibility of planing
the path after $d_t$ in $K$ raises.\\
For, during the first raise, we put negative the excursion
beginning at $t_0(w)$ and positive all the other excursions in $\left(t_0(w), t'\right)$
of height greater than $\varepsilon'$. Then during the second raise,
we put negative the excursion whose support contains this of
$de_{t_1}$, and so on. At the end of such $K$ raises, the path on
$[t_0, t']$ has an absolute value which doesn't exceed $\varepsilon'
+ K \varepsilon^{''}$, and $\varepsilon'$ on the excursion
straddling $t'$.

So we are going now to analyze its behavior on
$(t_d,t_{d+1})$. Let us denote $\gamma_{t_d}$ the first time after
$t_d$ at which one of the $\Gamma^\sigma$, $0\leq \sigma\leq s_0+l_0$,
vanishes on $(t_d,t_{d+1})$, and $\sigma_0$ the corresponding
level.\\

\begin{prop}
There exists a $\sigma(\Gamma^{\sigma_0^{'}+1})-$measurable,
$\mathbb{N}-$valued r.v., $K^{'}_{\sigma_0}$ such that there exists
$K^{'}_{\sigma_0} - 1$ r.v. $P_1, \dots, P_{K^{'}_{\sigma_0} - 1}$ themselves
with values in $\mathbb{N}$ and $\sigma (\Gamma^{\sigma_0 +
1})-$measurable, such that :
\begin{enumerate}
\item \label{l1} the $K^{'}_{\sigma_0} - 1$ excursion intervals
$e_{P_1}(\Gamma^{\sigma_0}), \dots, e_{P_{K^{'}_{\sigma_0} -
1}}(\Gamma^{\sigma_0})$ are disjoint and included in $\left(
d_{t_d}(\Gamma^{\sigma_0}), t_0(\Gamma^{\sigma_0}) \right)$
\item \label{l2} the heights $H_1, \dots, H_{K'-1}$ of these $K'-1$ excursions of
$\Gamma^{\sigma_0}$ satisfy on $[\sigma_0 \geq 0]$ :
\begin{equation*}
\left| \tau^0_{RL_{\sigma_0}-s_0 )}(\widetilde{w}) (\varphi(t_{d+1})) \right| - \varepsilon' < H_1 + \dots +
H_{K^{'}_{\sigma_0}-1} < \left|\tau^0_{RL_{\sigma_0}-s_0 )}(\widetilde{w}) (\varphi(t_{d+1}))\right| + \varepsilon'.
\end{equation*}
\end{enumerate}
\end{prop}
~\\
\begin{prf}
Noticing that the process :
\begin{equation*}
X = \sum_{n\in\mathbb{N}} \Gamma^n\mathds{1}_{[\sigma_0 = n]}
\end{equation*}
is absolutely continuous w.r.t. $\pi$, the proposition is an immediate consequence of lemma 5.
\end{prf}
~\\
Often in the sequel, we will denote $K^{\prime}_{\sigma_0}$ simply by K', if there is no ambiguity.\\
Set :
\begin{equation*}
\sigma'_0 = \sigma_0 - \left( K( \Gamma^{\sigma - 1}) + K'_{\sigma_0} \right).
\end{equation*}
\begin{prop}
For all $(n, \omega)$ such that $\sigma_0(\omega)
\leq n < \sigma_0(\omega) - K_{\sigma_0}(\omega) - K'_{\sigma_0}(\omega)$,
$\mathcal{N}^n(\omega)$ and $\Sigma^n(\omega)$ can be chosen so that
:
\begin{enumerate}
\item \label{l1} $0 \leq \Gamma^{\sigma_0 - K}_{t_{d+1}} <
\varepsilon'$ and $\|\widetilde{\Gamma}^{RL(\sigma_0 - 1) - s_0}_{|[0, t_d]}
- \Gamma^{\sigma_0}_{|[0, t_d]} \|_{\infty} < K \varepsilon''$.
\item \label{l2} $- (\mathbf{I} \circ \Gamma^{\sigma_0})(t_{d+1}) = H_{K'_{\sigma_0}-1}$
\item \label{l3} $\Gamma^{\sigma_0}$ has $K'-2$ tall excursions
included in $(t_d, t_{d+1})$ : $E_1 < E_2 < \dots < E_{K'-2}$ with
respective heights $H_1, \dots, H_{K'-2}$ verifying $\left|
i_{\Gamma^{\sigma_0}} E_1 \right| < \left| i_{\Gamma^{\sigma_0}} E_2 \right| < \dots < \left| i_{\Gamma^{\sigma_0}} E_{K'-2} \right|$.
\item \label{l4} $H_{n+1} + \dots + H_{K'_{\sigma_0}-1} \leq \Gamma^{\sigma_0 - K -n}_{t_{d+1}} < H_{n+1} + \dots + H_{K'_{\sigma_0}} +
\varepsilon' + n K' \varepsilon''$.
\item \label{l5} $H_{K'_{\sigma_0}-n} < -(\mathbf{I} \circ \Gamma^{\sigma_0 - K - n})(t_{d+1}) < H_{K'_{\sigma_0}-n} + K' \varepsilon''$.
\item \label{l6} $\Gamma^{\sigma_0 - K_{\sigma_0} -n}$ has $K'_{\sigma_0}-n-2$ tall excursions included in
$(t_d, t_{d+1}) E_1^n < \dots < E^n _{K'_{\sigma_0}-n-2}$ such that :
$\left|i_{\Gamma^{\sigma_0 - K_{\sigma_0} - n}} E_1 \right| < \dots < \left| i_{\Gamma^{\sigma_0 - K_{\sigma_0} - n} E_{K'_{\sigma_0}-n-2}} \right| <
\varepsilon''$, and whose heights $H^n_1, \dots, H_{K'_{J_0}-n-2}^n$
satisfy :
\begin{equation*}
H_l \leq H_l^n < H_l + n \varepsilon'' \ for \ n+1 \leq l < K'_{\sigma_0} - 1.
\end{equation*}
\end{enumerate}
\end{prop}
~\\

In our pursuit of the procedure, we can state :

\begin{prop}
It is possible to choose $\mathcal{N}^n(\omega)$ and
$\Sigma^n(\omega)$, for all $(n, \omega)$ such that $n = \sigma'_0$, in order to have :

\begin{equation*}
\left| \Gamma^{\sigma'_0}_{t_{d+1}} -
\tau^0_{RL(\sigma'_0)-s_0}(\varphi(t_{d+1})) \right| <
\varepsilon' + \frac{1}{2} (K'_{\sigma_0}-1) (K'_{\sigma_0}-2) \varepsilon''
\end{equation*}

\end{prop}
~\\
\begin{prf}
We notice that the building
excursions which appear in Proposition 5, the
$e_{P_0}(\Gamma^{\sigma_0 + K_{\sigma_0}})$'s, are successively protected. Once
protected, each of them receives a small excursion of height lower
than $\varepsilon''$ at each raise. So we deduce the result.
\end{prf}
~\\

\begin{prop}
It is possible to define $\mathcal{N}^n(\omega)$, $\Sigma^n(\omega)$
and $RL_{n}(\omega)$ inductively on the event $[n \geq \sigma'_0]$ if $\Gamma^n(\omega)$ doesn't vanish in $(t_d, t_{d+1})$ and
$RL_n(\omega) > s_0$, in such a way that :
\begin{enumerate}
\item \label{l1} $\left\|  \Gamma^n_{|[t_0(B^{\sigma_0}),t_{d+1}]} -
\Gamma^n_{t_{d+1}} \right\| _{\infty} < \varepsilon' + K(B^{\sigma_0})\varepsilon''$.
\item \label{l2} $\left| \Gamma^n_{t_{d+1}} -
\tau^0_{RL_n - s_0}(\varphi(t_{d+1})) \right| <
\varepsilon' + \frac{1}{2}(K^{\prime}_{\sigma_0} - 1)(K^{\prime}_{\sigma_0} - 2) \varepsilon^{''} + 2^{n - (\sigma_0-\sigma'_0)} \varepsilon''$.
\end{enumerate}
\end{prop}

~\\
\begin{prf}
We first recall that when we know $\Gamma^{n+1}(\omega)$, we know also
wether $\Gamma^n$ vanishes in $(t_d, t_{d+1})$. If it isn't the case, we
put :
\begin{equation*}
RL_n(\omega) = RL_{n+1}(\omega) - 1,
\end{equation*}
and define $\mathcal{N}^n(\omega)$ and $\Sigma^n(\omega)$ as we do
in previous propositions; but this time the counting of "errors" is
radically different. It can happen between $0$ and $t_d$ that an
excursion which was protected before become negative and, in the
following raise, is going to add to another protected excursion. So,
at each vertical raise, the "errors" are double of those of the
preceding raise.
\\
In the same manner the excursion straddling $t_{d+1}$ receives an
excursion with beginning in $[0, t_d]$. So, to the errors soon
acquainted at level $J_0^{''}$ we must add the error between $0$ and
$t_d$ of the preceding level which entails (ii).\\
For (i) : here the raises which are involved, are the planing one's, i.e. the
$K(\Gamma^{\sigma_0})$ first raises. At most, at each instant of the interval
$[t_0(B^{J_0 - 1}), t_{d+1}]$, the path $\Gamma^n$ has received $K(\Gamma^{\sigma_0})$ small
excursions. Then, this part of the path is just successively translated, which
entails (i).
\end{prf}
~\\
Then we define $\sigma_n, \ n \geq 0,\ hz_n \ and \
vt_n$ in the following manner :
\begin{equation*}
\forall n \in \mathbb{N}, \sigma_{n+1} = \sup \left(\sup \{ p > \sigma'_n, \Gamma^p \ vanishes \ in \ (t_d, t_{d+1}) \}, 0\right)
\end{equation*}
$hz_n$ (resp. $vt_n$) is the number of horizontal (resp. vertical) raisings occurring between levels $r$ and $n$.
As before, $S= \sup \{n \leq r; RL_n = s_0 \}$.

\begin{prop}
\begin{enumerate}
\item  For all $n \in \mathbb{N}$, $\sigma_n$, $hz_n$, $vt_n$ are r.v.
\item It is possible to define $\mathcal{N}^n$ and $\Sigma^n$ on the
event $\left[ \sigma_{k-1}^{'} \geq n \geq \sigma_k \right]$ in such a manner that
:
\begin{compactenum}
\item \label{l1} $\left\| \Gamma^n_{|[t_0 (B^{\sigma_{k-1}}), t_{d+1}]} -
\Gamma^n_{t_{d+1}} \right\|_{\infty} < \varepsilon' + K (\Gamma^{\sigma_k})
\varepsilon''$
\item \label{l2} $\left|\Gamma^n_{t_{d+1}} -
\tau^0_{RL_n-s_0}) (\varphi(t_{d+1})) \right| <
\varepsilon' + \frac{1}{2}(K^{\prime}_{\sigma_k} - 1)(K^{\prime}_{\sigma_k} - 2) \varepsilon^{''}  + hz_n \varepsilon'' 2^{vt_n}$
\end{compactenum}
\end{enumerate}
\end{prop}
~\\
\begin{prf}
It is the same as in the previous proposition.\\
For (a) we notice that, after the intervention of the planing excursions,
this part of the path is merely translated, without being affected by any other modification.
\end{prf}
~\\
\begin{prop}
\begin{enumerate}
\item \label{l1} $S$ is a r.v. such that $S - s_0$ doesn't depend upon $s_0$ and $\varepsilon''$, and we can choose $s_0$ large enough for $\mathbf{P}(A_2^{\varepsilon}) > 1 - \varepsilon(1+\frac{d}{d_0}) - \frac{\varepsilon}{d_0}$, where $A_2^{\varepsilon} := A_1^{\varepsilon} \cap [S \geq 0]$
\item \label{l2} Let $(\Gamma^n)$ be the sequence associated to the $\mathcal{N}^n$ and $\Sigma^n$.
It satisfies on $A^{\varepsilon}_2$ :

\begin{eqnarray*}
\left\| \Gamma^S_{|[0, t_d]} - \widetilde{\Gamma}^0_{|[0,t_d]}
\right\|_{\infty} &<& 2\varepsilon \\
\left\|  \Gamma^S_{|[t_0(B^{\sigma_{n_S}}), t_{d+1}]} -
\Gamma^S_{ t_{d+1}}
\right\|_{\infty} &<& \varepsilon' + K (B^{\sigma_{n_S}})\varepsilon'', \qquad where \ n_S = \sup \{n \leq r / \sigma_n \leq S \} \\
\left| \Gamma^S_{t_{d+1}} - \varphi(t_{d+1}) \right| &<&
2\beta + \varepsilon' + \frac{1}{2}(hz_n - 1)(hz_n - 2)
\varepsilon^{''} 2^{l_0} + \left(\varepsilon'' hz_{J_0 -
S}\varepsilon^{''} \right).
\end{eqnarray*}
\end{enumerate}

\end{prop}
~\\
\begin{prf}
For (i), see Lemma 12.\\
Then (ii) follows immediately from preceding Propositions.
\end{prf}
~\\
Then \textit{$\Gamma$ so defined is a $B-$raised Brownian motion.}\\
So, let us choose : $\varepsilon' = \frac{\varepsilon}{16}$
and $\varepsilon^{''} = \frac{\varepsilon}{16 (\frac{1}{2} (s_0 - 1)(s_0  - 2) + 2^{l_0} s_0  }$.\\
From the independence of the r.v. $S - s_0$, upon $s_0$ and
 $\varepsilon^{''}$, these choices don't create any vicious circle, and we can claim :

\begin{prop}
For all $\varepsilon > 0$, there exists a $B-$raised Brownian
motions verifying on the event $A^{\varepsilon}_2$ :
\begin{eqnarray}
\left\| \Gamma^0_{|[0, t_d]} - \varphi_{|[0, t_d]} \right\|_{\infty} &<&  \varepsilon \\
\left\| \Gamma^0_{|[t_0 (B^{\sigma_{n_S}}), t_{d+1}]} - \varphi_{|[t_0 (B^{\sigma_{n_S}}), t_{d+1}]} \right\|_{\infty} &<& \frac{\varepsilon}{8} + \frac{\varepsilon}{4} \\
\left| \Gamma^0_{t_{d+1}} - \varphi(t_{d+1}) \right| &<& \frac{\varepsilon}{4}.
\end{eqnarray}
\end{prop}
~\\
\begin{prf}
We deduce immediately these increases from Proposition 25 since $S
\leq s_0$, $hz_{J_0 - S} \leq s_0$, $K(\Gamma^{J_S-1}) \leq s_0$.
\end{prf}
~\\
At this point, the last task to achieve is to control $\Gamma^0$
between
 times $t_0 (B^{J_0})$ and $t_{d+1}$.\\

So we are going now to analyze more in details its behavior on
$(t_d,t_{d+1})$. Let us denote $\gamma_{t_d}$ the first time after
$t_d$ at which one of the $\Gamma^\sigma$, $0\leq \sigma\leq s_0+l_0$,
vanishes on $(t_d,t_{d+1})$, and $\sigma_0$ the corresponding
level.\\
Let us introduce the rectangle
$Rect_{\sigma_0}$ defined by the four straight lines with equations
:\\
$ x=t_d\;, x=t_{d+1}\;,
y=\inf\tau^0_{(RL_{\sigma_0} - s_0)} (\varphi)|_{[t_d,t_{d+1}]}
\;\,\,\qquad and\\
y=\sup\tau^0_{(RL_{\sigma_0} - s_0)} (\varphi)|_{[t_d,t_{d+1}]}
$
\\
$Rect_{\sigma_0}$ contains by definition the path of
$\tau^0_{(RL_{\sigma_0} - s_0)} (\varphi)|_{[t_d,t_{d+1}]}$
and, from the choice of $\alpha_0$, its height is lower than
$\dis\frac{\varepsilon}{4}$.\\
Now consider the path of $\Gamma(w)^{\sigma'_0}|_{[t_d,t_{d+1}]}$, it takes
one of the two forms given in the appendix.
\\
In the two cases by hypothesis, the total variation of $w^{s_0+
r_d}$ on $[t_d,t_{d+1}]$ is lower than $\dis\frac{\varepsilon}{2}$.
So, by lemma 5, and the definition of $\gamma_{t_d}$, it is greater
or equal to that of $\Gamma(w)^{\sigma_0}$, on $[t_d,\gamma_{t_d}]$.
Consequently the path $w^{\sigma_1}$ can move again from
$Rect_{\sigma_0}$ but at most from
$\frac{\varepsilon}{4} + \frac{\varepsilon}{2}$ on the same interval.\\
And rapidly, it is bound to join $\varphi$ in
$Rect_{\sigma_0}$ by the building
excursions, the flat part remaining flat.\\
Therefore, the rectangle $RR_{\sigma_0}$ with the same center and
vertical straight lines bordering it, and height that of
$Rect_{\sigma_0}+\frac{3\varepsilon}{2}$, contains the path
$w^{\sigma_1}|_{[t_d,t_{d+1}]}$.\\
During the following raises, the rectangle $Rect_{\sigma_0}$,
according to lemma 5, moves by isometry. We call $R_\sigma$ its new
positions, and likewise $RR_\sigma$ that of $RR_{\sigma_0}$.\\
We can easily check that, for all $\sigma < \sigma_1$ corresponding to
a vertical raise, the path $w^\sigma|_{[t_d,t_{d+1}]}$ is
contained in $RR_{\sigma}$.\\
Finally, for $w\in A^{\varepsilon}_2$ at
level $0$ we have the desired property :
\begin{equation*}
    \left\|\Gamma^0_{\nu}|_{[t_d,t_{d+1}]}-\varphi|_{[t_d,t_{d+1}]}\right\|_\infty < \varepsilon
\end{equation*}
 \\
So we have proved the following :

\begin{prop}
For all $\varepsilon > 0$, there exists a disjointed sum of $B-$raised Brownian motions such that :
 on $A^{\varepsilon}_3$,
\begin{eqnarray}
\left\| \Gamma^0 |_{[t_d,t_{d+1}]} - \varphi|_{[t_d,t_{d+1}]} \right\| &<& \varepsilon \\
\left| \Gamma^0_{t_{d+1}} - \varphi (t_{d+1}) \right| &<& \frac{\varepsilon}{4}.
\end{eqnarray}
\end{prop}

Thus we establish that $\mathcal{P}_{d+1}$ is true. So,
by induction, $\mathcal{P}_d$ is true for all $d \leq d_0$, and we can claim :

\begin{prop}
For all $\varepsilon > 0$, there exists a disjointed sum of $B-$raised Brownian motions such that :
\begin{equation*}
\mathbf{P} \left( \left[ \| \Gamma^0_{|[0, T]} - \varphi_{|[0, T]} \|_{\infty} < \varepsilon \right] \cap \left[
|\Gamma_T - \varphi(T) | < \frac{\varepsilon}{4} \right]\right) > 1 - 2 \varepsilon.
\end{equation*}
\end{prop}

Then we can apply Proposition 1 to $G = \left[ \| w_{|[0,T]} - \varphi_{|[0,T]} \|_{\infty} > \varepsilon \right]$.\\
So,
\begin{equation*}
\mathbf{P} \left( \forall n \geq 0, \left\| B^n_{|[0, T]} - \varphi_{|[0, T]} \right\|_{\infty} > \varepsilon \right) < 2 \varepsilon
\end{equation*}
We deduce immediately :
\begin{equation*}
\mathbf{P} \left( \forall n \geq 0, \left\| B^n_{|[0, T]} - \varphi_{|[0, T]} \right\|_{\infty} > \varepsilon \right) = 0.
\end{equation*}
But this property is true again when we replace $1$ by $a$, for all $a > 0$ :
This means :

\begin{theo}~\\
For almost every $\omega \in \Omega$, the orbit of $B(\omega)$ :
\begin{equation*}
    orb(B (\omega )) = \{B^n ( \omega ) ;\ , n \in \mathbb{N}\}
\end{equation*}
is dense in $W$, equipped with the topology of uniform convergence
on compact sets.
\end{theo}
~\\

Let us notice that if, in place of restrain ourselves
with the open
sets $B$, we have shown :\\
$\forall B$ closed set in $W$,
\begin{equation*}
\mathbf{P}(B)>0\Rightarrow \mathbf{P}(orb(w)\cap B\neq\emptyset)=1
\end{equation*}
Then every set $A$ $\mathbf{T}$-invariant, measurable and not negligible,
would contain the event $[orb(w)\cap B\,\neq\emptyset]$ and so,
would be
almost sure. Therefore, $\mathbf{T}$ would be ergodic.\\
To end, we are going to claim in an equivalent way, following thus
an interesting suggestion of J.P Thouvenot :
\begin{equation*}
    \forall(\varphi,\varepsilon)\in W|_{[0,1]}\times\mathbb{R}^+_\ast\,,
\end{equation*}
the reverse martingale $\mathbf{P}(w\in
B(\varphi,\varepsilon)|\mathcal{W}^n_\infty)$ admits a regular
conditional version $\mathbf{P}(w\in B(\varphi,\varepsilon)|w^n)$, and we
have :\\
\begin{theo}~\\
\begin{equation*}
    \mathbf{P}\;\textrm{a.s.}\,,\;\lim_{n\to\infty}\mathbf{P}(w\in
    B(\varphi,\varepsilon)|w^n)>0
\end{equation*}
\end{theo}
~\\
\noindent\begin{prf}\textbf{ of theorem 3.}\\
Suppose the contrary, and let :
\begin{equation*}
A:=\left[w\in W, \lim_{n\to\infty}\mathbf{P}(w\in
    B(\varphi,\varepsilon)|w^n)=0\right]
\end{equation*}
As $\mathbf{P}(w\in B(\varphi,\varepsilon)|w^n)=\mathbf{P}(w\in
 B(\varphi,\varepsilon)|w^{n+1})$, because $\mathbf{T}$ is measure-preserving. So we have :
\begin{equation*}
w\in A \Leftrightarrow \mathbf{T}w\in A
\end{equation*}
So $A$ is $\mathbf{T}$-invariant. Consequently :
\begin{equation*}
E\left(\mathds{1}_A \mathbf{P}(w\in
B(\varphi,\varepsilon)|w^n)\right)=\mathbf{P}\left(A\cap[w\in
B(\varphi,\varepsilon)]\right)=\mathbf{P}\left(A\cap[w^n\in
B(\varphi,\varepsilon)]\right)
\end{equation*}
But by hypothesis :
\begin{equation*}
\lim_{n\to\infty}E\left(\mathds{1}_A \mathbf{P} (w\in
B(\varphi,\varepsilon)|w^n)\right)=0
\end{equation*}
Therefore,
\begin{equation*}
\mathbf{P}\left(A\cap[orb(w)\cap
B(\varphi,\varepsilon)\,\neq\emptyset]\right)=0
\end{equation*}
which, from theorem 1, entails that $\mathbf{P}(A)=0$
\end{prf}~\\
\noindent Finally, let us remark that, if we could show :
\begin{equation*}
    \lim_{n\to\infty}\mathbf{P}\left(\big[w\in
    B(\varphi,\varepsilon)\big]|w^n\right)=\mathbf{P}\left(\big[w\in
    B(\varphi,\varepsilon)\big]\right)\;,
\end{equation*}
Than, not only $\mathbf{T}$ would be ergodic but exact which means :
\begin{equation*}
\mathcal{W}^\infty_\infty:=\underset{n\in\mathbb{N}}{\cap}\mathcal{W}^n_\infty
\textrm{ would be trivial.}
\end{equation*}

\end{document}